\documentclass{article}

\usepackage{amssymb,dsfont,latexsym}

\newtheorem{thm}{Theorem}

\newtheorem{lem}[thm]{Lemma}
\newtheorem{cor}[thm]{Corollary}

\newcommand\enu[1]{\smallskip\newline\makebox[5mm][l]{\rm(#1)}}

\newcommand\bp{\noindent{\it Proof.}\ }

\newcommand\1[1]{{\cal #1}}

\begin{document}

\author{Erling St{\o}rmer}

\date{9-12-2009 }

\title{Mapping cones of positive maps}

\maketitle
\begin{abstract}
We study mapping cones of positive maps  of $B(H)$ into itself, i.e. cones which
are closed under composition with completely positive maps.  As applications we
obtain characterizations of linear functionals with strong positivity properties with
respect to so-called symmetric mapping cones, with special emphasis on separable 
and PPT states.
\end{abstract}

\section*{Introduction}
A substantial part of the theory of positive maps of operator algebras is best understood
by considering cones of maps, especially the so-called mapping cones of positive maps of 
$B(H)$ into itself.  Given a mapping cone $\1 K$,   maps with strong positivity properties
with respect to $\1 K$, called $\1 K-$positive maps, are important and have been studied
in \cite{SSZ}\cite{S2}\cite{S3}\cite{S4}\cite{S5}.  In general it is not clear
which maps are $\1 K-$positive, but in \cite{S2}, Thm. 3.6, we showed a 
general result which implies that for a large class of mapping cones, called symmetric in the 
sequel, a map of $B(H)$ into itself is $\1 K-$positive if and
only if it belongs to $\1 K$.  In the present paper we shall give a more accessible proof of this latter result
via a proof which gives more insight into the theory of mapping cones than the previous 
proof.  Furthermore we shall show that the dual cone of a symmetric mapping cone is
itself a symetric mapping cone.  As applications of our results we get conditions for linear 
functionals to be positive on certain
cones of operators defined by mapping cones. We then comment on the connection with
separable states and obtain two new equivalent conditions for a state on 
$B(H)\otimes B(H) (= B(H\otimes H))$ to be a PPT-state. 

To be more specific let $H$ be a Hilbert space, which we for simplicity assume is finite
dimensional. Let $P(H)$ denote the set of positive linear maps of $B(H)$ into itself. 
By a \textit{mapping cone} $\1 K$ we mean a closed subcone of 
 $P(H)$   such that if $\phi\in \1 K,$ and $\alpha$ and $\beta\in
CP(H),$ the completely positive maps of $B(H)$ into itself, then $\alpha\circ\phi
 \in \1 K$ and $\phi\circ\beta\in \1K$.  We say $\1 K$ is \textit{symmetric} 
 if $\phi\in \1 K$ iff
$\phi^* \in \1 K$ iff $\phi^t = t\circ\phi\circ t\in \1 K,$ where $\phi^*$ is the
adjoint map of $\phi$ in the Hilbert-Schmidt structure on $B(H),$ i.e. $Tr(\phi(a)b)=
Tr(a\phi^*(b))$ for $a,b\in B(H),$ and $t$ is the transpose map on $B(H)$ with respect 
to an orthonormal base, and $Tr$ is the usual trace on $B(H).$ For a mapping cone $\1K$ 
we shall use the notation $\1K^*$ (resp $\1K^t$) for the cones consisting of $\phi^*$ 
(resp $\phi^t$) with $\phi\in\1K$, and we define the dual cone $\1K^\circ$ of $\1K$
to be the cone
$$
\1K^\circ = \{\phi\in P(H): Tr(C_{\phi}C_{\alpha})\geq 0\ \forall\ \alpha\in \1K\}.
$$
If $A$  is an \textit{operator system}, i.e. a closed self-adjoint linear subspace of $B(L)$ for a 
Hilbert space $L,$ we denote by 
$$
P(A,\1 K) = \{x\in A\otimes B(H): \iota\otimes\alpha (x)\geq 0 \ \forall\ \alpha\in\1 K\}.
$$
We say a map $\phi\colon A\to B(H)$ is \textit{$\1 K$-positive} if the linear functional
$\tilde\phi$ on $A\otimes B(H)$  defined by
$$
\tilde\phi(a\otimes b) = Tr(\phi(a)b^t)
$$
is positive on $P(A,\1K).$ By \cite{S2}, Lem.2.1, $\phi$ is positive iff $\tilde\phi$
is positive on the cone $A^+\otimes B(H)^+$ generated by tensors $a\otimes b$ with
$a\in A^+, b\in B(H)^+$, and by \cite{S2},Thm.3.2, $\phi$ is completely positive iff $\tilde\phi$ is a positive 
linear functional on $A\otimes B(H).$  

If the Hilbert space $L$ satisfies $dim L\leq dim H$, then we can consider $B(L)$ as
imbedded in $B(H)$ and thus, in order to study positive linear maps of $A$ into $B(H)$,
we can, and do, consider $A$ as a subset of $B(H)$.  We can now state our main results.

\begin{thm}\label{thm2}
Let $H$ be a finite dimensional Hilbert space and $\1K$ a symmetric mapping cone on
$B(H)$. Then its dual cone $\1K^\circ$ is also a symmetric mapping cone.
\end{thm}
\begin{thm}\label{thm1}
Let $H$ be a finite dimensional Hilbert space and $A\subset B(H)$ an operator system.
Let $\1 K$ be a symmetric mapping cone on $B(H),$ and let $\phi\colon A\to B(H)$
be a positive map.  Then $\phi$ is $\1 K-$positive if and only if $\phi$ is the restriction
to $A$ of a map in $\1 K$.
\end{thm}

Note that by \cite{S2},Thm.3.1, if $\phi$ is $\1 K-$positive then $\phi$ has a
$\1 K-$positive extension  $\psi\colon B(H)\to B(H)$.  Since the restriction to $A$
of a $\1 K-$positive map $\psi\colon B(H)\to B(H)$ is $\1 K-$positive, it suffices to
prove  Theorem 2 for $A = B(H).$  We shall therefore for the rest of the paper assume
$A=B(H)$ with $H$ finite dimensional.

\section {Proof of the theorems }
 If $(e_{ij})$ is a complete set of matrix units in $B(H)$ then the Choi matrix for a map $\phi$
is the operator
$$
C_{\phi} = \sum_{ij} e_{ij}\otimes \phi(e_{ij}) \in B(H\otimes H).
$$
By \cite{S3},Lem.5, $C_{\phi^t}$ is the density operator for $\tilde\phi,$ and by 
\cite{S5},Lem.4, $C_{\phi^t} = C_{\phi}^t.$  Furthermore, by \cite{Cho} $\phi$
is completely positive iff $C_{\phi}\geq 0,$ which holds iff  $C_{\phi}^t\geq 0.$  Thus 
$\tilde\phi \geq 0$ iff $\phi$ is completely positive, a result stated under more general 
situations in the introduction. 

\begin{lem}\label{lem3}
Let $dim H = n$ and $e_{1},\dots,e_{n}$ be an orthonormal basis for $H$. Let $J$ be 
the conjugation on $H\otimes H$ defined by 
$$
J z e_{i}\otimes e_{j}= \bar z e_{j}\otimes e_{i}
$$
with $z\in\0C$. Let $\phi\in P(H)$.  Then $C_{\phi^*}= JC_{\phi}J.$
\end{lem}
\bp
Let $V=(v_{ij})_{ij\leq n}\in B(H)$, and let $e_{ij}$ denote the matrix units 
corresponding to $e_{1},\dots,e_{n}$. Then  a straightforward computation yields
$$
AdV(e_{kl})= V^*e_{kl}V= (\bar v_{ki} v_{lj})_{ij}.
$$
Since $V^*=(\bar v_{ji})$ it follows that
$$
AdV^*(e_{kl})= Ve_{kl}V^* = (v_{ik}\bar v_{jl})_{ij}.
$$
From the definition of $J$ it thus follows that
\begin{eqnarray*}
JC_{AdV}J(z e_{p}\otimes e_{q}) &=& J(\sum_{ijkl} e_{kl}\otimes \bar v_{ki}v_{lj}e_{ij})(\bar z e_{q}\otimes e_{p})\\
&=&\sum_{ijkl}v_{ki}\bar v_{lj}e_{ij}e_{p}\otimes z e_{kl}e_{q}\\
&=&(\sum_{ijkl} v_{ik}\bar v_{jl} e_{kl} \otimes e_{ij})(z e_{p}\otimes e_{q})\\
&=& C_{AdV^*}(z e_{p}\otimes e_{q}),
\end{eqnarray*}
where we at the third equality sign exchanged (i,j) with (k,l). Since the vectors $e_{p}\otimes e_{q}$ 
form a basis for $H\otimes H$, $JC_{AdV}J=C_{AdV^*}.$  Now, if $\phi$ is a positive map, 
then $C_{\phi}$ is  a self-adjoint operator, hence the difference between two positive operators,
which both are the Choi matrices for completely positive maps.  Using the Kraus decompositions
for these completely positive maps, we see that $\phi$ is a real linear sum of maps $AdV$. Now the
adjoint map of $AdV$ is $AdV^*,$ as is easiy verified.  Applying this to each summand $AdV$, we 
thus get $JC_{\phi}J=C_{\phi^*}$, completing the proof of the  lemma.

\medskip

\textit{Proof of Theorem 1}

We first show $\1K^\circ =(\1K^{\circ})^t  $.  As remarked in the introduction $C_{\phi^t}=C_{\phi}^t.$
Thus if $\phi\in \1K^\circ, \alpha\in\1K,$ then
$$
0\leq Tr(C_{\phi}C_{\alpha})=Tr(C_{\phi^{tt}}C_{\alpha})=Tr(C_{\phi^t}^t C_{\alpha})=Tr(C_{\phi^t}C_{\alpha}^t)=Tr(C_{\phi^t}C_{\alpha^t}).
$$
Since $\1K$ is symmetric, $\1K=\1K^t.$  Thus $Tr(C_{\phi^t}C_{\alpha})\geq 0$ for all $\alpha\in\1K,$
hence $\phi^t \in \1K^\circ.$ The converse folows by symmetry of the argument.

In order to show $\1K^\circ = (\1K^{\circ})^*$ note that the map $\gamma (x)=Jx^*J$ is an
antiautomorphism of $B(H\otimes H),$ hence is trace invariant. In particular, if $x\in B(H\otimes H)$,
then
$$
Tr(JxJ) = Tr(\gamma(x^*)) = Tr(x^*).
$$
Thus by Lemma 3 if $\phi\in\1K^\circ, \alpha\in\1K$, then 
\begin{eqnarray*}
Tr(C_{\phi^*}C_{\alpha})&=& Tr(JC_{\phi}JC_{\alpha})\\
&=&Tr(JC_{\phi}JC_{\alpha}JJ)\\
&=&Tr((C_{\phi}C_{\alpha^*})^*),
\end{eqnarray*}
which is the complex conjugate of the positive number $Tr(C_{\phi}C_{\alpha^*})$, 
as $\alpha^*\in \1K$, so is itself positive. Thus $\phi^*\in \1K^\circ,$ completing the proof.
\medskip

We now embark on the proof of Theorem 2. For this we need to consider the cone $\1K^\sharp $
for $\1 K$  a mapping cone on $B(H).$  Following \cite{S5} we denote by $C^\1 K$ 
 the closed cone generated by all cones
 $$
 \iota\otimes\alpha^* (B(H\otimes H)^+), \ \alpha\in \1K.
 $$
 We denote by  $\1 K^\sharp$ the closed cone
 $$
   \1 K^\sharp = \{\beta\in P(H): \iota\otimes \beta (x)\geq 0 \ \forall\ x\in C^\1K \}.
  $$
Then $ \1 K^\sharp$
 is a mapping cone characterized by the property that a map $\phi$ belongs to 
 $ \1 K^\sharp$ iff $\phi\circ\alpha^*\in CP(H)$ for all $\alpha\in \1 K.$ 
 Furthermore, if $\1  K$ is symmetric then $C^{\1 K} = P(B(H;\1 K^\sharp).$ We denote by $\1P_{\1K}$
the cone of $\1K-$positive maps in $P(H).$  We are now ready to prove Theorem 1.  We
divide the proof into some lemmas.

\begin{lem}\label{lem2}
Let $\1K$ be a symmetric mapping cone on $B(H).$  Then we have
\enu{i} $\1K \subset \1P_{\1K}$.
\enu{ii} $\1 K^{\sharp} = (\1 K^{\sharp} )^t.$
\enu{iii} If  $\1 K^{\sharp} = (\1 K^{\sharp} )^*$ then $\1 K^{\sharp} = \1P_{\1K^{\sharp}}$.
\end{lem} 
\bp
Let $\pi\colon B(H)\otimes B(H) \to B(H)$ by $\pi(a\otimes b) = b^t a.$  Then
by \cite{S5}, Lem.10, $Tr\circ \pi$ is positive, and
$$
\tilde\phi = Tr\circ\pi\circ (\iota\circ\phi^{*t}).
$$
Since $\1K$ is symmetric, if $\phi\in \1K,$ then $\phi^{*t}\in \1K$, so that $\tilde\phi$
is positive on $P(B(H),\1K),$ i.e. $\phi$ is $\1K-$positive, proving (i).

Ad(ii).  We first note that if $\alpha,\beta \in B(B(H),B(H))$ then 
$$
(\alpha\circ\beta)^t = \alpha^t \circ\beta^t.
$$ 
Indeed, if $x\in B(H),$ then
$$
(\alpha\circ\beta)^t (x) =( (\alpha\circ\beta)(x^t))^t = \alpha(\beta(x^t))^t =\alpha^t(\beta(x^t)^t)= \alpha^t (\beta^t (x)),
$$
proving the assertion.

To show (ii) let $\phi\in P(H).$ By the characterization of $\1 K^{\sharp}$ mentioned
before the statement of the lemma, $\phi\in \1 K^{\sharp}$ iff $\phi\circ \alpha^*\in CP(H)$
for all $\alpha\in\1K$ iff $\phi\circ \alpha\in CP(H)$ for all $\alpha\in \1K,$ since 
$\1K =\1K^*.$  Since a map $\beta\in P(H)$ belongs to $CP(H)$ iff $\beta^t \in CP(H)$,
it follows from the above assertion that $\phi\in\1 K^{\sharp} $ iff $\phi^t\circ\alpha^t =
(\phi\circ\alpha)^t \in CP(H)$ for all $\alpha\in\1K,$ and since $\1K = \1K^t,$ iff 
$\phi^t \circ\beta \in CP(H)$ for all $\beta\in \1K,$ i.e. iff $\phi^t\in \1K^{\sharp},$ completing 
the proof of (ii). 

Ad(iii). Assume further that $\1K^{\sharp}=(\1K^{\sharp})^*.$  By (ii) $\1K^{\sharp}=(\1K^{\sharp})^t$,
so $\1K^{\sharp}$ is a symmetric mapping cone.  Therefore by (i) applied to $\1K^{\sharp}$, 
$\1K^{\sharp} \subset \1P_{\1K^{\sharp}}.$  By \cite{S5},Thm.12, $\1P_{\1K^{\sharp}}=
\1P_{\1K}^{\circ}.$  Let $\phi\in \1P_{\1K}^{\circ}.$  
By \cite{S5},Thm.1, $\alpha\circ\phi\in CP(H) $ for all $\alpha\inÊ\1K^t = \1K.$
Hence $\phi^*Ê\circ\alpha^* = (\alpha\circ\phi)^* \in CP(H)$ for all $\alpha\in \1K.$ 
Thus $\phi^*\in\1K^{\sharp}=(\1K^{\sharp})^*,$ so that $\phi\in\1K^{\sharp}$. Thus 
 $\1P_{\1K^{\sharp}}=\1P_{\1K}^{\circ}\subset\1K^{\sharp},$ completing the proof of
 the lemma.
 
 \medskip
 It should be remarked that the proof of \cite{S5},Thm.1, (i)$\Leftrightarrow$(ii) makes use of
 \cite{S2},Thm.3.6.  But that is unnecessary in our setting as seen from the proof of Lemma 7
 below.
 
 The proof of the next lemma does not use finite dimensionality, so we state it for general $C^*-$algebras.
 We denote by $\1P_{\1K}(M)$ the set of $\1K-$positive maps of a $C^*-$algebra $M$ into $B(H).$
 
 \begin{lem}\label{lem3}
 Let $M$ be a $C^*-$algebra and $\1K\subset\1L$ be mapping cones on $B(H)$ such that
  $\1P_{\1K}(M) =   \1P_{\1L}(M).$ Then $P(M,\1K) = P(M,\1L)$.
  \end{lem}
 \bp
 From the definition of $P(M,\1K)$ it is clear that  $\1K\subset\1L$ implies  $P(M,\1K) \supset  P(M,\1L)$.
 If $\phi\colon B(H)\to B(H)$, and $\tilde\phi$ is positive on $P(M,\1L)$, then $\phi\in \1P_{\1L}(M) =   \1P_{\1K}(M),$
 so $\tilde\phi$ is positive on $P(M,\1K)$. Since the identity operator $1$ is an interior point of $ P(M,\1L),$
 the Hahn-Banach Theorem for cones implies that  $P(M,\1K) = P(M,\1L)$, proving the lemma.
 
 \begin{lem}\label{lem4}
 Let $\1K$ be a symmetric mapping cone. Then
 \enu{i} $\1K\subset (\1K^{\sharp})^{\sharp}.$
 \enu{ii}$ P(B(H),\1K)=P(B(H), (\1K^{\sharp})^{\sharp}).$
 \end{lem}
 \bp
 We have
 $$
 \1K^{\sharp} = \{\beta\in P(H): \beta\circ\alpha^*\in CP(H) \ \ \forall \alpha\in \1K\}.
 $$ 
Let $\phi\in \1K$.  Then $\beta\circ\phi^* \in CP(H)$ for all $\beta\in \1K^{\sharp},$ hence
$\phi\circ\beta^* = (\beta\circ\phi^*)^* \in CP(H)$ for all $\beta\in \1K^{\sharp}.$  Thus by
the above formula for   $\1K^{\sharp}$ applied to $(\1K^{\sharp})^{\sharp}, \phi\in (\1K^{\sharp})^{\sharp}$,
proving (i).  

(ii) follows from from Lemma 5, since by \cite{S5},Thm.12
$\1P_{\1K}= (\1P_{\1K^{\sharp}})^\circ =\1P_{(\1K^{\sharp})^{\sharp}},$ using that by Lemma 4
$\1K^{\sharp} = (\1K^{\sharp})^t.$  The proof is complete.

\begin{lem}\label{lem5}
Let $\1K$ and $\1L$ be mapping cones on $B(H)$ such that  $\1K\subset\1L$ and 
$\1P_{\1K} = \1P_{\1L}$.  Then $\1K = \1L.$
\end{lem}
\bp
By Lemma 5  $P(B(H),\1K) = P(B(H),\1L)$.  Each operator in   $P(B(H),\1K)$ is of the form
$C_{\phi}$ with $\phi\colon B(H)\to B(H).$  By definition of  $P(B(H),\1K)$ we then have, since $\1K^*$
is a mapping cone,
\begin{eqnarray*}
C_{\phi}\in P(B(H),\1K)&\Leftrightarrow& \iota\otimes\alpha(C_{\phi})\geq 0\  Ê\forall\ \alpha \in \1K\\
&\Leftrightarrow&Tr((\iota\otimes\alpha)(C_{\phi})C_{\psi})\geq 0\   \forall\ \alpha\in \1K, \psi\in CP(H)\\
 &\Leftrightarrow&Tr(C_{\phi}C_{\alpha^*\circ\psi})= Tr(C_{\phi} (\iota\otimes\alpha^*)(C_{\psi}))\geq 0\ \forall\ \alpha\in \1K, \psi\in CP(H)\\
  &\Leftrightarrow&Tr(C_{\phi}C_{\beta})\geq 0\ \forall\ \beta\in \1K^*\\
  &\Leftrightarrow&\phi\in (\1K^*)^\circ.
\end{eqnarray*}
Since  $P(B(H),\1K) = P(B(H),\1L)$, and the same equivalence as the one above holds for $\1L,$ we
have $(\1K^*)^\circ =(\1L^*)^\circ $.  Thus $\1K^* =(\1K^*)^{\circ\circ}=(\1L^*)^{\circ\circ}=
\1L^*.$  But the map $\phi\to \phi^*$ is a bijection of $\1K$ onto $\1K^*,$  so $\1K = \1L,$ 
completing the proof.

\begin{lem}\label{lem6}
Let $\1K$ be a symmetric mapping cone on $B(H).$  Then $\1K =(\1K^{\sharp})^{\sharp}. $
\end{lem}
\bp
This is immediate from Lemmas 6 and 7

\medskip
\textit{Completion of the proof of Theorem 2}

As remarked after the statement of Theorem 2 it suffices to prove it for $A=B(H).$ 
By Lemma 8, since $\1K$ is symmetric, so that $(\1K^{\sharp})^{\sharp}$ is
symmetric, by using Lemma 4, \cite{S5},Thm.1, and Thm.12 in that order, we have
\begin{eqnarray*}
\1K &=& \1K^{\sharp})^{\sharp}\\
&=&\{\beta\in P(H): \beta^*\circ\alpha^*\in CP(H) \  \forall\ \alpha\in \1K^{\sharp}\}Ê\\
&=&\{\beta\in P(H): \alpha\circ\beta \in CP(H)  \  \forall\ \alpha\in \1K^{\sharp}=(\1K^{\sharp})^t\}\\
&=&(\1P_{ \1K^{\sharp}})^\circ\\
&=&\1P_{ \1K}.
\end{eqnarray*}
The proof is complete.
\medskip

It thus follows from Theorem 1 that $\1P_{\1K}$ and $\1P_{\1K}^\circ $ are symmetric mapping cones when $\1K$
is.  The next corollary shows the same for $\1K^\sharp.$

\begin{cor}\label{cor1}
Let $\1K$ be a symmetric mapping cone on $B(H).$ Then $\1K^{\sharp}=\1K^\circ.$  
\end{cor}
\bp
By \cite{S5}, Lem.7 and Thm.1, and then Theorem 1 we get
\begin{eqnarray*}
\1K^{\sharp}&=& \{\beta\in P(H): \beta\circ\alpha^*\in CP(H) \ \forall\ \alpha\in \1K\}\\
&=&\{\beta\in P(H): \alpha\circ\beta^*\in CP(H) \ \forall\ \alpha\in \1K\}\\
&=&\{\beta^*\in P(H): \alpha\circ\beta\in CP(H) \ \forall\ \alpha\in \1K\}\\
&=&\{\beta^*\in P(H): \beta\in (\1P_{\1K})^\circ \}\\
&=&\{\beta^*\in P(H): \beta\in \1K^\circ\}\\
&=&(\1K^\circ)^*\\
&=&\1K^\circ.
\end{eqnarray*}
The proof is complete.

\section{Linear functionals}
In this section we apply the theorems to linear functionals.  The first result is closely related to
\cite{S5},Theorem 1, applied to the dual cone $\1K^\circ$ of $\1K.$  Then we consider 
applications to separable and PPT-states on $B(H)\otimes B(H) (= B(H\otimes H).$

\begin{thm}\label{thm2}
Let $\1K$ be a symmetric mapping cone on $B(H),$ and let $\rho$ be a linear functional on
$B(H\otimes H)$ with density operator $h$. Then the following conditions are equivalent.
\enu{i} $\rho = \tilde\phi$ with $\phi\in \1K^\circ.$
\enu{ii} $\rho(C_{\alpha})\geq 0$ for all $\alpha\inÊ\1K.$
\enu{iii} $\iota\otimes\alpha(h)\geq 0$ for all $\alpha\inÊ\1K.$
\enu{iv} $\rho\circ (\iota\otimes\alpha)\geq 0$ for all $\alpha\inÊ\1K.$
\enu{v} $\rho$ is positive on the cone $\{x\in B(H\otimes H): \iota\otimes\alpha(x)\geq 0\  \forall\ \alpha\inÊ\1K^\circ\}$.
\end{thm}
\bp
By \cite{S2},Lem.2.1 $\rho =  \tilde\phi$ with $\phi\colon B(H)\to B(H).$ By
\cite{S3}, Lem.5, and \cite{S5}, Lem.4, we have $h = C_{\phi^t} = C_{\phi}^t.$  Let
$\alpha\in \1K,$ then we have
$$
\rho(C_{\alpha}) = Tr(h C_{\alpha}) = Tr( C_{\phi}^t C_{\alpha}).
$$
Since $\1K = \1K^t$ it follows that (ii) holds iff $\phi\in \1K^\circ,$ and by Theorem 2
iff $\phi\in \1P_{\1K}^\circ.$  An application of \cite{S5}, Thm.1, now shows the
equivalence of (i) with (iii) and (iv).  To show the equivalence of (iv) and (v) note that (iv) holds iff 
$\rho$ is positive on the cone generated by $\iota\otimes\alpha (B(H\otimes H)^+)$ for all
$\alpha\in \1K$, i.e. $\rho$ is positive on $C^\1K$.  But by \cite{S5}, Lem.11, 
$$
C^\1K = P(B(H),\1K^\sharp) = \{x\in B(H\otimes H): \iota\otimes\alpha(x)\geq 0\  \forall\ \alpha\inÊ\1K^\sharp\}.
$$
Thus an application of Corollary 9 completes the proof.

\medskip
\textit{Remark 9} By \cite{S5},Thm.1(iv), it follows that the dual cone of a mapping cone is itself a mapping cone.
In particular this holds for $P(H) $ and the cone $S(H)$ of superpositive, also called entanglement
breaking maps, generated by maps $\phi(x)=\omega(x)a,$ where $\omega$ is a state on $B(H)$
and $a\in B(H)^+.$  By \cite{S2}, Lem.2.4, $S(H)$ is the minmal and $P(H)$ the maximal
mapping cones on $B(H),$ hence $S(H) = P(H)^\circ$.  We have by \cite{HSR}, or by \cite{S3},
Thm.2, that a state $\rho = \tilde\phi$ is separable iff $\phi\in S(H) = P(H)^\circ$.  Thus by Theorem 10
applied to $\1K = P(H)$ we see that a state $\rho$ on $B(H\otimes H) $ is separable iff it is positive on the 
cone $C$  consisting of $x\in B(H\otimes H)$ such that $\iota\otimes\omega(x)\geq 0$ for all states  $\omega$  on  $B(H)$.
 This result is true for states on $A \otimes B(H)$ with $A$ an operator system, see \cite
{S4},Prop.1, and has recently been extended to operator spaces in \cite{VTT}. Note that 
$x = \sum a_{i}\otimes b_{i}\in B(H\otimes H)$ belongs to the cone $C$ above iff $\eta(x)\geq 0$
for all separable states $\eta$ of $B(H\otimes H)$.  Indeed, letting $\omega$ and $\eta$ denote states on $B(H)$ we have
\begin{eqnarray*}
x\in C &\Leftrightarrow& 0\leq \iota\otimes\omega( \sum a_{i}\otimes b_{i}) = \sum a_{i}\omega(b_{i})\otimes 1 \ \forall\ \omega \\
&\Leftrightarrow& \sum a_{i}\omega(b_{i}) \geq 0  \ \forall\ \omega \\
&\Leftrightarrow& \sum \eta(a_{i})\omega(b_{i}) = \eta(\sum a_{i}\omega(b_{i}))\geq 0 \ \forall\ \eta, \omegaÊ\\
&\Leftrightarrow& \eta\otimes\omega (x)\geq 0 \ \forall \eta, \omega \\
&\Leftrightarrow& \tau(x)\geq 0 
\end{eqnarray*}
for all separable states $\tau $ on $B(H\otimes H)$.
 
The reader should also notice that condition (iii) in Theorem 10 corresponds to the Horodecki Theorem \cite{HHH}
for separable states.

\medskip
If $A$ is an operator system a state $\rho$ on $A\otimes B(H) $ is called a \textit{PPT-state} if
$\rho\circ(\iota\otimes t)$ is also a state.  Let $\1P$ denote the symmetric mapping cone
consisting of maps $\phi$ in $P(H)$ which are both completely positive and copositive, where the
latter means that $t\circ\phi$ is completely positive.  It is known \cite{S3},Prop.4, that $\rho$
is PPT iff $\rho = \tilde\phi$ with $\phi\in \1P$, and if $A= B(L)$ for some Hilbert space $L,$ by 
\cite{S3},Thm.8, iff $\iota\otimes t (h)\geq 0,$  where $h$ is the density operator for $\rho.$
Using Theorem 2 we can add two more equivalent conditions for a state on $B(H)\otimes B(H)$
to be PPT.
\begin{thm}\label{thm3}
Let $H$ be a finite dimensional Hilbert space and $\rho$ a linear functional on $B(H)\otimes B(H)$.
Then the following conditions are equivalent.
\enu{i} $\rho$ is a PPT-state.
\enu{ii} $\rho$ is positive on the set $E=\{x\in B(H\otimes H): x\geq 0$ or $\iota\otimes t(x)\geq 0\}$.
\enu{iii}  $\rho$ is positive on the cone $\{x\in B(H\otimes H):\iota\otimes \alpha (x)\geq 0 \ \forall\ \alpha\in \1P\}$.

\end{thm}
\bp
As mentioned before $\rho=\tilde\phi$ with $\phi\colon B(H)\to B(H)$, and $\rho$ is a PPT-state
iff $\phi\in\1P,$ which by \cite{S5},Lem.14, holds iff $\rho\geq 0$ on $E.$ Thus $(i)\Leftrightarrow
(ii).$  Condition (iii) means that $\phi$ is $\1P-$positive, which by Theorem 2 means that $\phi\in\1P,$
which again is equivalent to $\rho$ being a PPT-state.  This proves $(i)\Leftrightarrow (iii).$  The proof is
complete.

\medskip
A variation of condition (ii) in Theorem 11 is \cite{S4},Cor.19, where it is shown that $\phi$ is
decomposable, or equivalently $\rho=\tilde\phi$ is PPT, iff $\rho$ is positive on the cone
$\{x\in B(H\otimes H)^+: \iota\otimes t(x)\geq 0\}.$

Department of Mathematics, University of Oslo, 0316 Oslo, Norway.

e-mail: erlings@math.uio.no

\end{document}